\documentclass[11pt]{article}

\usepackage{amsmath}
\usepackage{amsthm}
\usepackage{amscd}
\usepackage{amssymb}
\usepackage{verbatim}

\title{Real-analytic, volume-preserving actions of lattices on
$4$-manifolds}
\author{Benson Farb and Peter Shalen \thanks{Both authors are supported 
in part by the NSF.}}

\theoremstyle{plain}
\newtheorem{theorem}{Theorem}
\newtheorem{proposition}[theorem]{Proposition}

\newtheorem{conjecture}[theorem]{Conjecture}
\newtheorem{question}[theorem]{Question}

\def\proof{{\bf {\medskip}{\noindent}Proof. }}

\def\endproof{$\diamond$ \bigskip}

\def\title{\em}

\newcommand\R{\mbox{\bf R}}

\newcommand\Z{\mbox{\bf Z}}
\newcommand\Q{\mbox{\bf Q}}

\newcommand\Diff{\rm Diff}

\newcommand\Rrank{\R\mbox{-rank}}
\newcommand\Qrank{\Q\mbox{-rank}}

\DeclareMathOperator\SL{SL}
\DeclareMathOperator\GL{GL}

\begin{document}
\maketitle
\begin{abstract}
We prove that any real-analytic, volume-preserving action of a lattice
$\Gamma$ in a simple Lie group with $\Qrank(\Gamma)\geq 7$ on a closed
$4$-manifold of nonzero Euler characteristic factors through a finite
group action.  
\end{abstract}

\section{Results}

Zimmer conjectured in \cite{Zi1} that the standard action of $\SL(n,Z)$ on 
the $n$-torus is minimal in the following sense:

\begin{conjecture}[Zimmer]
\label{conjecture:zimmer}
Any smooth, volume-preserving action of 
any finite-index subgroup $\Gamma<\SL(n,\Z)$  on a 
closed $r$-manifold factors through a finite group action if $n>r$.
\end{conjecture}

While Conjecture \ref{conjecture:zimmer} has been proved for actions 
which also preserve an extra geometric structure such as a
pseudo-Riemannian metric (see, e.g. \cite{Zi1}),  
almost nothing is known in the general case.  For $r=2$ and $n>4$, the
conjecture was proved for real-analytic actions in \cite{FS1}.  Quite 
recently, Polterovich \cite{Po} has brought ideas from symplectic
topology to the problem, using these to give a proof of 
Conjecture \ref{conjecture:zimmer} for
orientable surfaces of genus $>1$.  In \cite{FS2} we will point out how
his methods actually prove Conjecture \ref{conjecture:zimmer} for the
torus as well. For $r=3$, Conjecture \ref{conjecture:zimmer} is known 
only in some special cases where $\Gamma$ contains some torsion 
and the action is real-analytic (see
\cite{FS1}). 

In this note we prove the following result.

\begin{theorem}[Actions on $4$-manifolds with \boldmath$\chi(M) \neq 0$]
\label{theorem:main}
Let $\Gamma$ be a lattice in a simple Lie group $G$ such that 
$\Qrank(\Gamma)\geq 7$.  Then any 
real-analytic, volume-preserving action of $\Gamma$ on a 
closed $4$-manifold of nonzero Euler characteristic 
factors through a finite group action.  
\end{theorem}

In particular, Theorem \ref{theorem:main} applies to any finite-index
subgroup of $\SL(n,\Z)$, $n>7$.  

The main ingredient in the proof of Theorem 
\ref{theorem:main} is Theorem 7.1 of \cite{FS1} on real-analytic actions
which preserve a volume form.  This theorem, which is the most
difficult result in \cite{FS1}, gives a {\em codimension-two} invariant
submanifold for centralizers of elements with fixed-points.  One can 
then apply results of \cite{FS1} and \cite{Re}, which show that 
real-analytic (not necessarily area-preserving) 
actions of certain lattices on $2$-dimensional manifolds must factor
through finite groups.

For the case of symplectic actions, some further progress on Conjecture 
\ref{conjecture:zimmer} can be found in \cite{Po} and \cite{FS2}.

\section{Proof of Theorem \ref{theorem:main}}

Before giving the proof of Theorem \ref{theorem:main}, we will need two 
algebraic properties of lattices with large $\Qrank$.

\begin{proposition}[Some algebraic properties of lattices]
\label{proposition:alg}
Let $\Gamma$ be a lattice in a simple algebraic group over $\Q$.  
Then the following hold:
\begin{enumerate}
\item If $d=\Qrank(\Gamma)\geq 7$ then $\Gamma$ 
contains commuting subgroups $A$ and $B$
which are isomorphic to irreducible lattices with $\Qrank$s $2$ and
$d-3$ respectively.
\item If $\Qrank(\Gamma)\geq 4$ then $\Gamma$ contains a torsion-free 
nilpotent subgroup which is not metabelian.
\end{enumerate}
\end{proposition}

\proof  The proof of the first statement 
is similar to that of Proposition 2.1 of \cite{FS1}. By Margulis's
Arithmeticity Theorem (see, e.g., \cite{Zi2}, Theorem 6.1.2), $\Gamma$
is commensurate with the group of $\Z$-points of a simple algebraic
group $G$ defined over $\Q$.  Hence without loss of generality we
can assume that $\Gamma$ itself is the group of $\Z$-points in such a
group $G$.    

Since $G$ is simple, the root system $\Phi$ of $G$ is irreducible, and
the Dynkin diagram determined $\Phi$ therefore appears in the list given
in Section 11.4 of [Hu]. By going through this list, one sees that in
every case where $d\ge7$, one may ``erase a vertex $v$'' of the diagram
to obtain a a graph with 2 components: one with two vertices and another
which is a Dynkin diagram with at least $d-3$ vertices.  Let $G_1$ and
$G_2$ be the root subgroups corresponding to these two components of the
Dynkin diagram.  Then the group of $\bf Q$-points of $G_1$ 
has $\Qrank$ at least $2$, the group of $\Q$-points of $G_2$ has
$\Qrank$ at least $d-3$, and $G_1$ commutes with $G_2$.

Now $\Gamma_i=\Gamma\cap G_i$ is an 
arithmetic lattice in $G_i$ for $i=1,2$, since by a theorem of 
Borel-Harish-Chandra (see, e.g. \cite{Zi2}) the $\Z$-points of an 
algebraic group defined over 
$\Q$ form a lattice in the group of $\R$-points.  Then $A=\Gamma_1$ and 
$B=\Gamma_2$ have the required properties.

To prove the second statement, note that since $\Qrank(\Gamma)\geq 4$,
we can find a connected, nilpotent Lie subgroup $N$ which is defined
over $\Q$ and has derived length $3$, i.e. is not metabelian. As
$\Gamma\cap N$ is the group of $\Z$-points of the $\Q$-group
$N$, it is a lattice in $N$, and in particular is
Zariski-dense in $N$. Hence $\Gamma\cap N$ is nilpotent and has no
metabelian subgroup of finite index. As $\Gamma\cap N$ must have a
tosrion-free subgroup of finite index, the assertion follows.
\endproof

We now turn to the proof of Theorem \ref{theorem:main}.  Let $M$ be a closed
$4$-manifold with nonzero Euler characteristic.  Let $\Gamma$ be an
irreducible lattice in a simple Lie group $G$, and assume
$d=\Qrank(\Gamma)\geq 7$. By part (1) of Proposition
\ref{proposition:alg}, $\Gamma$ contains commuting subgroups $A$ and $B$
which are isomorphic to irreducible lattices with $\Qrank$s $2$ and
$d-3\geq 4$ respectively.

Let $\gamma_0$ be any infinite order element of $A$.  By an old theorem
of Fuller \cite{Fu}, any homeomorphism of a closed manifold of nonzero
Euler characteristic has a periodic point; the proof is an application
of the Lefschetz fixed-point theorem and basic number theory.  Hence
some positive power $\gamma$ of $\gamma_0$ has a fixed point.  

We will also need the following two facts.  First, since $\Qrank(B)\geq
d-3\geq 4$, it follows from Margulis's Superrigidity Theorem that any
representation of $B$ into $\GL(4,\R)$ has finite image.  Second, since
$\Gamma$ is a lattice in a simple Lie group $G$ with $\Rrank(G)\geq 2$,
the Margulis Finiteness Theorem (see, e.g., Theorem 8.1 of \cite{Zi2})
gives that $\Gamma$ is {\em almost simple} in the sense that 
any normal subgroup
of $\Gamma$ must be finite or of finite index.

We are now in a position to apply Theorem 7.1 of \cite{FS1}.  For the
reader's convenience we recall the statement here.  We say that a group
action $\rho: \Gamma\to \Diff(M)$ is {\em infinite} if $\rho$ has
infinite image.

\bigskip
\noindent
{\bf Theorem 7.1 of \cite{FS1}: }
{\it Let $\Gamma$ be an almost simple group.  Suppose we are given an
infinite, volume-preserving, real-analytic action of $\Gamma$ on a
closed, connected $n$-manifold $M$.  Suppose further that $\Gamma$
contains commuting subgroups $A$ and $B$ with the following properties:
\begin{itemize}
\item There exists an element $\gamma\in A$, noncentral in $\Gamma$, 
having a fixed point in $M$.
\item $A$ is isomorphic to an irreducible lattice of $\Q$-rank $\geq 2$.
\item $B$ is noncentral in $\Gamma$.
\item Any representation of any finite-index subgroup of 
$B$ in $\GL(n,\R)$ has finite image.
\end{itemize}

Then there is a nonempty, connected, 
real-analytic submanifold $W\subset M$ of 
codimension at least 2 which is invariant under a finite-index subgroup 
$B'$ of $B$.  Furthermore, the action of this subgroup on $W$ is infinite.
}

\bigskip
\noindent
{\bf Remark. }The action of $B'$ on the 
surface $W$ produced by this theorem is NOT necessarily 
area preserving.  
\bigskip

We now conclude the proof of Theorem \ref{theorem:main}.  Since
$B'$ is a lattice in a simple Lie group and 
$\Qrank(B')\geq 4$, it follows from part (2) of Proposition
\ref{proposition:alg} that $B'$ contains 
a torsion-free nilpotent subgroup $H$ which is not metabelian.  But
Rebelo \cite{Re} showed that any nilpotent group of real-analytic
diffeomorphisms of a compact, oriented surface must be metabelian.  It
follows that the action of $H$ on $W$ is not effective. 

Since $H$ is
torsion-free, there is an infinite-order element of $H\le B'$ which acts
trivially on $W$. Since $B'$ has finite index in the almost simple group
$B$, and hence is almost simple, 
some finite index subgroup $C$ of $B'$ acts trivially on
$W$; in particular $C$ has a global fixed point in $M$.  Since 
$C$ is isomorphic to a lattice of $\Qrank$ at least $4$, by 
Lemma 3.2 of \cite{FS1} we have that 
a finite index subgroup $D$ of $C$ acts trivially on $M$.  Since
$\Gamma$ is almost simple, it follows that some finite index
normal subgroup of $\Gamma$ acts trivially on $M$, and we are done.
\endproof

\noindent
Benson Farb:\\
Dept. of Mathematics, University of Chicago\\
5734 University Ave.\\
Chicago, Il 60637\\
E-mail: farb@math.uchicago.edu
\medskip

\noindent
Peter Shalen:\\
Dept. of Mathematics, University of Illinois at Chicago\\
Chicago, IL 60680\\
E-mail: shalen@math.uic.edu

\end{document}